\documentclass[12pt]{amsart}

\usepackage{amsmath,amsthm,amssymb}

\font\teneufm=eufm10 \font\seveneufm=eufm7 \font\fiveeufm=eufm5
\newfam\frakturfam
\textfont\frakturfam=\teneufm \scriptfont\frakturfam=\seveneufm
\scriptscriptfont\frakturfam=\fiveeufm
\def\frak{\fam\frakturfam}

\newtheorem{lm}{Lemma}
\newtheorem{theor}{Theorem}
\newtheorem{co}{Corollary}

\def\bee{\begin{eqnarray}}
\def\bes{\begin{eqnarray*}}
\def\eee{\end{eqnarray}}
\def\ees{\end{eqnarray*}}

\def\a{\alpha}
\def\b{\beta}
\def\g{\gamma}

\def\Proof{{\sl Proof.}\ }

\newcommand{\Aut}{{\mathrm{Aut}}}
\newcommand{\ad}{\mbox{ad}}

\newcommand{\id}{\mbox{id}}
\newcommand{\de}{\partial}
\newcommand{\J}{\mbox{J}}

\title{Automorphisms of elliptic Poisson algebras}

\begin{document}
\date{}
\maketitle

\begin{center}
{\textbf{Leonid Makar-Limanov} \footnote{Supported
by an NSA grant and by grant FAPESP, processo  06/59114-1.}}\\
Department of Mathematics \& Computer Science,\\
Bar-Ilan University, 52900 Ramat-Gan, Israel and\\
Department of Mathematics, Wayne State University, \\
Detroit, MI 48202, USA\\
e-mail: {\em lml@math.wayne.edu}\\

and\\

{\bf Umut Turusbekova, \ \ Ualbai Umirbaev}\\
Department of Mathematics, Eurasian National University\\
 Astana, 010008, Kazakhstan \\
e-mail: {\em umut.math@mail.ru, \ \ umirbaev@yahoo.com}

\end{center}

\begin{abstract}
We describe the automorphism groups of elliptic Poisson algebras
on polynomial algebras in three variables and give an explicit set
of generators and defining relations for this group.
\end{abstract}


\noindent


\section{Introduction}

\hspace*{\parindent}

 It is well known \cite{Czer,Jung,Kulk,Makar} that the automorphisms
of polynomial algebras and free associative algebras in two
variables are tame. It is also known \cite{Umi25,Um19} that
polynomial algebras and free associative algebras in three
variables in the case of characteristic zero have wild
automorphisms. It was recently proved \cite{MLTU} that the
automorphisms of free Poisson algebras in two variables over a
field of characteristic $0$ are tame. Note that the Nagata
automorphism \cite{Nagata,Umi25} gives an example of a wild
automorphism of a free Poisson algebra in three variables.

One of the main problems of affine algebraic geometry (see, for
example \cite{Essen}) is a description of automorphism groups of
polynomial algebras in $n\geq 3$ variables. Unfortunately this
problem is still open even when $n=3$ and there isn't any
plausible conjecture about the generators. There was a conjecture
that the group of automorphisms of polynomial algebras in three
variables is generated by all affine and exponential automorphisms
(see \cite{Essen}). This seems to be not true and the first author
and D.\,Wright (oral communication) independently constructed
potential  counterexamples.

In order to find a plausible conjecture it is necessary to
consider many different types of automorphisms. In this paper we
describe the groups of automorphisms of the polynomial algebra
$K[x,y,z]$ over a field $K$ of characteristic $0$ endowed with
additional structure, namely, with Poisson brackets. A study of
automorphisms of Poisson structures on polynomial algebras is also
interesting in view of M.\,Kontsevich's Conjecture about the
existence of an isomorphism between automorphism groups of
symplectic algebras and Weyl algebras \cite{Belov-Kontsevich}.

A complete description of quadratic Poisson brackets on the
polynomial algebra \\
$K[x,y,z]$ over a field $K$ of characteristic $0$ is given in
\cite{Don}, \cite{DH}, and \cite{LX}. Among corresponding Poisson
algebras the most interesting are elliptic Poisson algebras
$E_{\a}$. By definition (see, for example \cite{Odesskii}), the
elliptic Poisson algebra $E_{\a}$  is the polynomial algebra
$K[x,y,z]$ endowed with the Poisson bracket defined by \bes
\{x,y\}=-\a xy+z^2, \ \ \{y,z\}=-\a yz +x^2, \ \ \{z,x\}=-\a
zx+y^2, \ees where $\a \in K$.

We describe the automorphism groups of the elliptic Poisson
algebras $E_{\a}$ over a field $K$ of characteristic $0$. We also
show that $E_{\a}$ doesn't have any nonzero locally nilpotent
derivations.

This paper is organized as follows. In Section 2 we give some
definitions and examples. In Section 3 we describe Casimir
elements of elliptic algebras and prove a lemma related to the
non-rationality of elliptic curves. In Section 4 we study  the
automorphism group of $E_{\a}$.

\section{Definitions and examples}

\hspace*{\parindent}

A vector space $P$ over a field $K$ endowed with two bilinear
operations $x\cdot y$ (a multiplication) and $\{x,y\}$ (a Poisson
bracket) is called {\em a Poisson algebra} if $P$ is a commutative
associative algebra under $x\cdot y$, $P$ is a Lie algebra under
$\{x,y\}$, and $P$ satisfies the following identity: \bes \{x,
y\cdot z\}=\{x,y\}\cdot z + y\cdot \{x,z\}. \ees

Let us call a linear map $\phi:P \longrightarrow P$ an automorphism of
$P$ as a Poisson algebra if
\bes \phi(xy)=\phi(x)\phi(y), \ \ \
\phi(\{x,y\})=\{\phi(x),\phi(y)\} \ees for all $x, y\in P$

Similarly, a linear map $D:P \longrightarrow P$ is a derivation of
$P$ if \bes D(xy)=D(x)y+xD(y), \ \ \
D(\{x,y\})=\{D(x),y\}+\{x,D(y)\} \ees for all $x, y\in P$.
 In other words, $D$ is simultaneously a
derivation of $P$ as an associative algebra and as a Lie algebra.
It follows that for every $x\in P$ the map
\bes
\ad_x : P \longrightarrow
P, \ \ (y\mapsto \{x,y\}),
\ees
is a derivation
of $P$. It is natural to call these
derivations \textit{inner}.

There are two important classes of Poisson algebras.

1) Symplectic algebras $S_n$. For each $n$ algebra $S_n$ is a polynomial algebra \\
$K[x_1,y_1, \ldots,x_n,y_n]$ endowed with the Poisson bracket
defined by \bes \{x_i,y_j\}=\delta_{ij}, \ \ \{x_i,x_j\}=0, \ \
\{y_i,y_j\}=0, \ees where $\delta_{ij}$ is the Kronecker symbol
and $1\leq i,j\leq n$.

2) Symmetric Poisson algebras $PS(\frak{g})$. Let ${\frak{g}}$ be
a Lie algebra with a linear basis $e_1,e_2,\ldots,e_k,\ldots$.
Then $PS(\frak g)$ is the usual polynomial algebra $K[e_1,e_2,\ldots,e_k,\ldots]$
endowed with the Poisson bracket defined
by
\bes \{e_i,e_j\}=[e_i,e_j]
\ees
 for all $i,j$, where $[x,y]$
is the multiplication of the Lie algebra $\frak{g}$.

Let $K\{x_1,x_2,\ldots,x_n\}$ be a free Poisson algebra in the
variables $x_1,x_2,\ldots,x_n$. Recall that
$K\{x_1,x_2,\ldots,x_n\}$ is the symmetric Poisson algebra
$PS(\frak g)$, where $\frak g$ is the free Lie algebra in the
variables $x_1,x_2,\ldots,x_n$. Choose a linear basis \bes e_1,
e_2, \ldots, e_m, \ldots \ees of $\frak{g}$ such that $ e_1=x_1,
e_2=x_2,\ldots, e_n=x_n$. The algebra $K\{ x_1,x_2,\ldots,x_n\}$
is generated by these elements as a polynomial algebra.
Consequently, \bes K[x_1,x_2,\ldots,x_n]\subset
K\{x_1,x_2,\ldots,x_n\}. \ees

This inclusion was successfully used in \cite{Umi24} to define
algebraic dependence of two elements of the polynomial algebra
$K[x_1,x_2,\ldots,x_n]$. Namely,  two elements $f$ and $g$ of
$K[x_1,x_2,\ldots,x_n]$ are algebraically dependent if and only if
$\{f,g\}=0$. It is also proved \cite{MakarU2} that if elements $f$
and $g$ of the free Poisson algebra $K\{x_1,x_2,\ldots,x_n\}$
satisfies the equation $\{f,g\}=0$ then $f$ and $g$ are
algebraically dependent. There is a conjecture \cite{MakarU2} that
if  $\{f,g\}\neq 0$ then the Poisson subalgebra of
$K\{x_1,x_2,\ldots,x_n\}$ generated by $f$ and $g$ is a free
Poisson algebra in the variables $f$ and $g$, i.e., the elements
$f$ and $g$ are {\em free}.

\section{Elliptic Poisson algebras}

\hspace*{\parindent}

Let \bes C=C(x,y,z)=\frac{1}{3} (x^3+y^3+z^3)-\a xyz, \ \ \ \ \a
\in K, \ees then the equation $C(x,y,z)=0$ defines an elliptic
curve ${\frak E}_{\a}$ in $KP^2$. The elliptic Poisson algebra
$E_{\a}$ (see, for example \cite{Odesskii}) is the polynomial
algebra $K[x,y,z]$ endowed with the Poisson bracket \bes
\{x,y\}=\frac{\partial C}{\partial z}, \ \ \{y,z\}=\frac{\partial
C}{\partial x}, \ \ \{z,x\}=\frac{\partial C}{\partial y}. \ees
Consequently, \bee\label{f2} \{x,y\}=-\a xy+z^2, \ \ \{y,z\}=-\a
yz +x^2, \ \ \{z,x\}=-\a zx+y^2. \eee
This bracket can be written
as $\{u, v\} = \J(u, v, C(x,y,z))$ where $\J(u, v, C(x,y,z))$ is
the Jacobian, i. e. the determinant of the corresponding Jacobi
matrix.

If $\a^3=1$ and $K$ contains a root $\epsilon$ of the equation $\lambda^2+\lambda+1=0$ then
\bee\label{f3}
C=C(x,y,z)=\frac{1}{3}(x+y+\a z)(\epsilon x+\epsilon^2 y+\a z)(\epsilon^2 x+\epsilon y+\a z),
\eee
It is well known (see, for example \cite{Huse}) that $C$ is an irreducible polynomial if $\a^3\neq 1$.

The center  $Z(A)$ of any Poisson algebra $A$ can be defined as a
set of all elements $f\in A$ such that $\{f,g\}=0$ for every $g\in
A$.  The elements of $Z(A)$ are called {\em Casimir} elements. The
statement of the following lemma is well known (see, for example
\cite{Odesskii}).

\begin{lm}\label{l5}
$Z(E_{\a})= K[C]$.
\end{lm}
\Proof Since $\{C, h\} = \J(C, h, C) = 0$ for any $h \in E_{\a}$
we see that $C\in Z(E_{\a})$. Assume that $Z(E_{\a}) \ni f$ which
is algebraically independent with $C$. Take an element $g \in
E_{\a}$ which is algebraically independent over $K[C, f]$. Then
for any element $h \in E_{\a}$ there exists a polynomial
dependence $H(C, f, g, h) = 0$. So $0 = \{g, H(C, f, g, h)\} =
\{g, h\} \frac{\de\, H}{\de\, h}$ and $\{g, h\} = 0$ since we can
assume that $H(C, f, g, h)$ has minimal possible degree in $h$.
Therefore $g$ is also in the center. Take now an element $p \in
E_{\a}$. Then $0 = \{p, H(C, f, g, h)\} = \{p, h\}
\frac{\de\,H}{\de\,h}$. Therefore $\{p, h\} = 0$ and $E_{\a}$ is a
commutative Poisson algebra, which is not the case. So  if $f \in
Z(E_{\a})$ then $f$ is algebraically dependent with $C$. Since the
brackets are homogeneous we can assume that $f$ is homogeneous.
But then $f = kC^n$ where $k \in K$. Indeed, $f$ and $C$ satisfy a
polynomial relation $P(C, f) = 0$ which can be assumed
homogeneous. So over an algebraic closure of $K$ we will have $f^m
= k C^n$ where $n$ and $m$ are relatively prime integers and $k
\in K$ since $f, \, C \in K[x, y, z]$. If $\a^3 \neq 1$ then $C$
is irreducible and so $m = 1$; if $\a^3 = 1$ then $C$ is either
irreducible or a product of three linear factors and again $m =
1$. $\Box$

\begin{lm}\label{l4}
Suppose $a, \, b, \, c \in K[x_1,x_2,\ldots,x_n]$ are homogeneous
elements of the same degree. If $C(a,b,c)=0$ and $\a^3\neq 1$ then
$a,\, b$, and $c$ are proportional to each other.
\end{lm}
\Proof The statement of the lemma is an easy corollary of
non-rationality of elliptic curves (see, for example \cite{Huse}).
But we prefer to give here an independent proof for the
convenience of the reader. So, \bee\label{f4} a^3+b^3+c^3-3\a
abc=0. \eee If $a$ and $b$ are divisible by an irreducible
polynomial $p$ then $c$ is also divisible by $p$ and we can cancel
$p^3$ out of (\ref{f4}). So we can assume that $a,\, b$, and $c$
are pair-wise relatively prime. If $\deg(a) \neq 0$ we can also
assume that $\frac{\partial a}{\partial x_1}\neq 0$. For any $f$
we put $f'=\frac{\partial f}{\partial x_1}$. Differentiating
(\ref{f4}), we get \bes a'a^2+b'b^2+c'c^2-\a (a'bc+ab'c+abc')=0.
\ees Consequently, \bes
a'(a^3+b^3+c^3-3\a abc)-a [a'a^2+b'b^2+c'c^2-\a (a'bc+ab'c+abc')]\\
= (b^2-\a ac)(a'b-ab')+(c^2-\a ab)(a'c-ac')=0, \ees and
\bee\label{f5} (b^2-\a ac)(a'b-ab')=(c^2-\a ab)(ac'-a'c). \eee
Note that $b^2-\a ac\neq 0$ and $c^2-\a ab \neq 0$, otherwise $a$,
$b$, and $c$ are not pair-wise relatively prime. If $b^2-\a ac$
and $c^2-\a ab$ are relatively prime then (4) gives that $b^2-\a
ac$ divides $ac'-a'c$ . But $\deg(ac'-a'c) < \deg(b^2-\a ac)$,
hence $ac'-a'c = 0$. Then $(a/c)' = 0$ and $a = cr$ where $r \in
K(x_2, \dots, x_n)$. But then $a$ and $c$ again are not relatively
prime: any irreducible factor of $a$ which contains $x_1$ should
divide $c$.

Consequently, $b^2-\a ac$ and $c^2-\a ab$ are not relatively prime. Let $p$ be an irreducible
polynomial that divides $b^2-\a ac$ and $c^2-\a ab$. Denote by $I$ the ideal of $K[x_1,x_2,\ldots,x_n]$
generated by $p$. For $f\in k[x_1,x_2,\ldots,x_n]$ denote by $\bar{f}$
the homomorphic image of $f$ in $K[x_1,x_2,\ldots,x_n]/I$.
Then, $\bar{b}^2= \a \bar{a}\bar{c}$ and $\bar{c}^2= \a \bar{a}\bar{b}$.
Since $a$, $b$, and $c$ are pair-wise relatively prime, $\bar{a}, \bar{b}, \bar{c}\neq 0$.
Note that $K[x_1,x_2,\ldots,x_n]/I$ is a domain.
So
$\bar{a}^2= \a \bar{b}\bar{c}$ and $\bar{b}^2 \bar{c}^2 = \a^2 \bar{a}^2 \bar{b}\bar{c} = \a^3 \bar{b}^2\bar{c}^2.$
This gives $\a^3=1$ contrary to our assumptions. $\Box$

\section{Automorphisms}

\hspace*{\parindent}

Let us denote by $\varphi_{\gamma}$, where $\gamma\in K^*$, an automorphism of $E_{\a}$ such that
\bes
\varphi_{\g}(x)=\g x, \ \ \varphi_{\g}(y)=\g y, \ \ \varphi_{\g}(z)=\g z.
\ees
Note that $\langle \varphi_{\g} \rangle \cong K^*$, where $K^*$ is the multiplicative group of the field $K$. Algebra
$E_{\a}$ has also automorphisms
\bes
\tau : E_{\a}\longrightarrow E_{\a}, \ \ \ x\mapsto y, \ \ y\mapsto z, \ \ z\mapsto x.
\ees
and

\bes
\sigma : E_{\a}\longrightarrow E_{\a}, \ \ \ x\mapsto x, \ \ y\mapsto \epsilon y, \ \ z\mapsto \epsilon^2 z,
\ees
if a solution $\epsilon$ of the equation $\lambda^2+\lambda+1=0$ is in $K$.

Denote by $\Aut\,E_{\a}$ the group of automorphisms of $E_{\a}$.
Since the center of $E_{\a}$ is $K[C]$ the restriction of $\psi$
to $K[C]$ is an automorphism of $K[C]$ and so $\psi(C) = k_1C +
k_2$ where $k_1, \, k_2 \in K$ and $k_1 \neq 0$. Let $G$ be the
group of those  automorphisms of $K[x, y, z]$ which preserve
$K[C]$.

\begin{lm}\label{l6} The group $G$  consists of linear automorphisms.
\end{lm}
\Proof Let $\psi \in G$. Then by the chain rule \bes \J(\psi(x),
\psi(y), \psi(C)) = \J(\psi(x), \psi(y), \psi(z)) \frac{\de \,
\psi(C)}{\de \, \psi(z)}  = k_1(\psi(z^2) - \a \psi(x) \psi(y)).
\ees Also \bes \J(\psi(x), \psi(y), \psi(C)) = \J(\psi(x),
\psi(y), k_2C + k_3) = k_2 \{ \psi(x), \psi(y) \}. \ees

So $k_2\{\psi(x), \psi(y)\} =  k_1(\psi(z^2) - \a\psi(x)\psi(y))$
where $k_1, \, k_2 \in K$ are non-zero constants. Similarly,
$k_2\{\psi(y), \psi(z) \} = k_1(\psi (x^2) - \a \psi (y) \psi(z))$
and $k_2\{\psi(z), \psi(x)\} = k_1(\psi(y^2) - \a\psi(z)\psi(x))$.
Denote $\psi(x)=a,\, \psi(y)= b,\, \psi(z)=c$. So \bee\label{f6}
k_4\{a,b\}=-\a ab+c^2, \ \ k_4\{b,c\}=-\a bc +a^2, \ \
k_4\{c,a\}=-\a ca+b^2. \eee for some non-zero $k_4 \in K$.

Using automorphism $\tau$ if necessary, without loss of generality
we can assume that $\deg(a) \leq \deg(c)$ and $\deg(b) \leq
\deg(c)$. Then (\ref{f6}) show that $\deg(a) = \deg(b) = \deg(c)$.
Denote by $\bar{f}$ the highest homogeneous part of $f\in
K[x,y,z]$. \bee\label{f7} k_4\{\bar{a},\bar{b}\}=-\a
\bar{a}\bar{b}+\bar{c}^2, \ \ k_4\{\bar{b},\bar{c}\}=-\a
\bar{b}\bar{c} +\bar{a}^2, \ \ k_4\{\bar{c},\bar{a}\}=-\a
\bar{c}\bar{a}+\bar{b}^2, \eee as the terms of the highest
possible degree.

If $\deg\,c>1$ then $C(\overline{a}, \overline{b}, \overline{c}) =
0$ since otherwise $\deg\,C(a, b, c) > 3$. So $\overline{a}, \,
\overline{b}$, and $\overline{c}$ are pair-wise proportional by
Lemma \ref{l4}. Therefore \bes
\{\bar{a},\bar{b}\}=\{\bar{b},\bar{c}\}=\{\bar{c},\bar{a}\}=0.
\ees Then (\ref{f7}) gives \bes \a \bar{a}\bar{b}=\bar{c}^2, \ \
\a \bar{b}\bar{c} =\bar{a}^2, \ \ \a \bar{c}\bar{a}=\bar{b}^2 \ees
and \bes \bar{b}^2 \bar{c}^2 = \a^2 \bar{a}^2 \bar{b}\bar{c} =
\a^3 \bar{b}^2\bar{c}^2. \ees This is impossible if $\a^3\neq 1$.
Therefore $\deg(a) = \deg(b) = \deg(c) = 1$.

If $\a^3=1$ then using the decomposition (\ref{f3}) we have $3 =
\deg(\psi(C)) = \deg(a + b+ \a c) + \deg(\epsilon a +\epsilon^2
b+\a c) + \deg(\epsilon^2 a+\epsilon b+\a c)$. Therefore $\deg(a +
b+ \a c) = \deg(\epsilon a +\epsilon^2 b+\a c) = \deg(\epsilon^2
a+\epsilon b+\a c) = 1$ and the elements $a$, $b$, and $c$ also
have degree $1$.

So, we proved that $\psi$ is an affine automorphism. If $\a^3
\not= 1$ then (\ref{f6}) gives the linearity of $\psi$ (just look
at the lowest homogeneous parts of $a, \,b$, and $c$). If $\a^3 =
1$ then non of the $\epsilon^{i}a +\epsilon^{2i} b+\a c$ can
contain a constant term since otherwise $\psi(C)$ contains a term
of degree less than three. So $\psi$ is a linear automorphism in
this case as well. $\Box$

\begin{theor}\label{t1}
Let $K$ be an arbitrary field of characteristic $0$ such that the
equation $\lambda^2 + \lambda + 1 = 0$ has a solution in $K$. If
$\a^3 \not= 1$ then the group of automorphisms $\Aut\,E_{\a}$ of
the algebra $E_{\a}$ is generated by $\varphi_{\g} (\g\in K^*)$,
$\sigma$, and $\tau$.
\end{theor}
\Proof
Let $\psi$ be an arbitrary automorphism of $E_{\a}$.
By Lemma \ref{l6}, $\psi$ is a linear automorphism. Let
\bes
\psi(x)=a=\b_{11} x+\b_{12} y +\b_{13} z, \\
\psi(y)=b=\b_{21} x+\b_{22} y +\b_{23} z, \\
\psi(z)=c=\b_{31} x+\b_{32} y +\b_{33} z.
\ees
Denote by $B$ the matrix $(\b_{ij})_{3\times 3}$.
Since
\bee\label{f8}
\{a,b\}=-\a ab+c^2, \ \ \{b,c\}=-\a bc +a^2, \ \
\{c,a\}=-\a ca+b^2.
\eee
we have
\bes
\{a,b\}=(\b_{11}\b_{22}-\b_{12}\b_{21}) \{x,y\}+(\b_{12}\b_{23}-\b_{13}\b_{22}) \{y,z\}+(\b_{13}\b_{21}-\b_{11}\b_{23}) \{z,x\}\\
=-\a (\b_{11}\b_{22}-\b_{12}\b_{21}) xy+(\b_{11}\b_{22}-\b_{12}\b_{21}) z^2 -\a (\b_{12}\b_{23}-\b_{13}\b_{22}) yz\\
+ (\b_{12}\b_{23}-\b_{13}\b_{22}) x^2 -\a (\b_{13}\b_{21}-\b_{11}\b_{23}) xz + (\b_{13}\b_{21}-\b_{11}\b_{23}) y^2.
\ees
Also
\bes
-\a ab+c^2= (\b_{31}^2-\a \b_{11}\b_{21}) x^2+(\b_{32}^2-\a \b_{12}\b_{22}) y^2\\
+(\b_{33}^2-\a \b_{13}\b_{23}) z^2 +[2\b_{31}\b_{32}-\a (\b_{11}\b_{22}+\b_{12}\b_{21})] xy\\
+[2\b_{32}\b_{33}-\a (\b_{12}\b_{23}+\b_{13}\b_{22})] yz
+[2\b_{31}\b_{33}-\a (\b_{11}\b_{23}+\b_{13}\b_{21})] zx. \ees
Using this and (\ref{f8}) we can write 18 relations between the
coefficients of the matrix $B$: \bee\label{f9} \b_{ij}^2=\a
\b_{i+1j}\b_{i+2j}+(\b_{i+1j+1}\b_{i+2j+2}-\b_{i+1j+2}\b_{i+2j+1})
\eee and \bee\label{f10} \b_{ij}\b_{ij+1}=\a \b_{i+1j+1}\b_{i+2j}
\eee for all $i$ and $j$ modulo $3$. Denote by $b^j$ the product
of all elements of the $j$th column of $B$. Then (\ref{f10}) gives
\bes b^jb^{j+1} = \a^3 b^{j+1}b^j \ees
So $b^jb^{j+1} = 0$ since
$\a^3\neq 1$ and at least two columns contain zero elements.
Therefore the matrix $B$ contains zero entries.

If $\a = 0$ then $\b_{ij}\b_{ij+1}= 0$ and every row contains two
zero entries. Of course, a row should contain a non-zero element
since the determinant of $B$ is not zero.

Assume now that $\a \neq 0$. If in some row we have three non-zero
elements we can check using (\ref{f10}) that all entries of $B$
are not equal to zero. Indeed, if $\b_{ij}\b_{ij+1} \neq 0$ then
$\b_{i+1j+1}\b_{i+2j} \neq 0$ and using this inequality for $j =
1, \, 2, \, 3$ we fill all the matrix with non-zero entries. So
each row has at least one zero element.

Suppose that every row contains just one zero. Up to renumbering
we can assume that $\beta_{11}\beta_{12} \neq 0$ and $\beta_{13} =
0$. Then $\beta_{22}\beta_{31} \neq 0$. Now, $\beta_{22}\beta_{23}
= \a \beta_{33}\beta_{12}$ and $\beta_{21}\beta_{22} = \a
\beta_{32}\beta_{11}$. If $\beta_{23} = 0$ then $\beta_{33} = 0$
which is impossible since $B$ has a non-zero determinant. Hence
$\beta_{23} \neq 0$, $\beta_{33} \not= 0$ and $\beta_{21} = 0$,
$\beta_{32} = 0$.  Then $\beta_{11}\beta_{12}= \a
\beta_{22}\beta_{31}$, $\b_{22}\b_{23}=\a \b_{33}\b_{12}$, and
$\b_{33}\b_{31}=\a \b_{11}\b_{23}$ in force of (\ref{f10}). If we
multiply these equalities we will get $\a^3 = 1$. So there exist a
row with two zeros.

If there is a row with two non-zero entries we can assume that
$\beta_{11}\beta_{12} \neq 0$. Then $\beta_{22}\beta_{31} \neq 0$
and it is easy to check that a non-zero element in any other
position will lead to two non-zero entries in every row. But then
we have a column of zeros which is impossible. So every row (and
every column) contains just one non-zero element.

Let us assume that $\b_{11} \neq 0$. Then the relations (\ref{f9})
give $\b_{11}^2 = (\b_{22}\b_{33}-\b_{23}\b_{32})$ and either
$\b_{22}\b_{33} = 0$ or $\b_{23}\b_{32} = 0$. If $\b_{11}^2 =
-\b_{23}\b_{32}$, then $\b_{23}^2 = -\b_{32}\b_{11}$ and
$\b_{32}^2 = -\b_{11}\b_{23}$ and we obtain a contradiction by
multiplying these equalities. So $\b_{11}^2 = \b_{22}\b_{33}$,
$\b_{22}^2 = \b_{33}\b_{11}$, and $\b_{33}^2 = \b_{11}\b_{22}$. We
have three solutions to these equations: $\b_{11} =
\b_{22}=\b_{33}$, $\b_{22}=\epsilon\b_{11},
\b_{33}=\epsilon^2\b_{11}$, and $\b_{22}=\epsilon^2\b_{11},
\b_{33}=\epsilon\b_{11}$. They correspond to the cases $\psi=
\varphi_{\b_{11}}$, $\psi=\sigma\varphi_{\b_{11}}$, and
$\psi=\sigma^2\varphi_{\b_{11}}$, respectively. Remaining two
solutions with $\b_{12}\b_{23}\b_{31} \neq 0$ and
$\b_{13}\b_{23}\b_{31} \neq 0$ correspond to additional actions by
$\tau$. So, this gives the statement of the theorem. $\Box$

Recall that a derivation $D$ of an algebra $R$ is called
\textit{locally nilpotent} if for every $a\in R$ there exists a
natural number $m=m(a)$ such that $D^m(a)=0$. If $R$ is generated
by a finite set of elements $a_1,a_2,\ldots,a_k$ then it is well
known (see, for example \cite{Essen}) that a derivation $D$ is
locally nilpotent if and only if there exist positive integers
$m_i$, where $1\leq i\leq k$, such that  $D^{m_i}(a_i)=0$.

\begin{co}\label{c2}
Algebra $E_{\a}$ does not have any nonzero locally nilpotent derivations.
\end{co}
\Proof Let $D$ be a nonzero locally nilpotent derivation of
$E_{\a}$. For every $c\in Z(E_{\a})$ and $f\in E_{\a}$ we have
\bes 0=D\{c,f\}=\{D(c),f\}+\{c,D(f)\}=\{D(c),f\}. \ees
So $D(c)\in
Z(E_{\a})$, i.e., $Z(E_{\a})$ is invariant under the action of
$D$. Therefore $D$ induces a locally nilpotent derivation of
$Z(E_{\a}) = K[C]$ and $D(C)=\b \in K$ (see, for example
\cite{Freu}). Since \bes D(C)=D(x) \frac{\partial C}{\partial
x}+D(y) \frac{\partial C}{\partial y}+D(z) \frac{\partial
C}{\partial z}=\b \ees and $C$ is homogeneous of degree three,
$\b=0$.

So $D(C)=0$ and $CD$ is also a locally nilpotent derivation. For
any locally nilpotent derivation its exponent is an automorphism
(see, for example \cite{Freu}). So $\exp(CD)$ is an automorphism
which is certainly not linear. This contradicts Lemma \ref{l6}.
$\Box$

Usually it is much easier to describe the locally nilpotent
derivations than the automorphisms. But at the moment we do not
know any direct proof of Corollary \ref{c2}.

By Theorem \ref{t1}, the automorphism group $\Aut\,E_{\a}$ of the
algebra $E_{\a}$ is generated by $\varphi_{\g} (\g\in K^*)$,
$\sigma$, and $\tau$. Note that $\langle \varphi_{\g}
\rangle=\{\varphi_{\g} | \g \in K^*\} \cong K^*$,  $\langle \tau
\rangle\cong Z_3$, and $\langle \sigma \rangle\cong Z_3$, where
$Z_3$ is the cyclic group of order $3$.   The automorphisms
$\varphi_{\g}$, where $\g\in K^*$, are related by \bee\label{f11}
\varphi_{\g_1} \varphi_{\g_2}=\varphi_{\g_1\g_2}. \eee Moreover,
it is easy to check that \bee\label{f12} \sigma \varphi_{\g}=
\varphi_{\g}\sigma, \tau \varphi_{\g}=\varphi_{\g} \tau, \ \ \
\eee i.e., $\varphi_{\g}$ belongs to the center of $\Aut\,E_{\a}$.

In addition,
\bee\label{f13}
\sigma^3=\id, \tau^3=\id, \sigma \tau=\varphi_\epsilon \tau \sigma.
\eee

\begin{theor}\label{t2}
Let $K$ be an arbitrary field of characteristic $0$ such that the
equation $\lambda^2+\lambda+1=0$ has a solution in $K$. Then, \bes
\Aut\,E_{\a}=(\langle \varphi_{\g} \rangle\times \langle \tau
\rangle) \rtimes \langle \sigma \rangle =(\langle \varphi_{\g}
\rangle\times \langle \sigma \rangle) \rtimes \langle \tau \rangle
\cong (K^*\times Z_3)\rtimes Z_3. \ees if $\a^3 \neq 1$
\end{theor}
\Proof Using (\ref{f11}), (\ref{f12}), and (\ref{f13}), any
automorphism $\psi\in \Aut\,E_{\a}$ can be written as \bes
\psi=\varphi_{\g}\tau^i\sigma^j \ees
where $\g \in K^*$ and $0\leq
i,j\leq 2$. It is easy to check that this representation is
unique. The relations (\ref{f11})-(\ref{f13}) show that $\langle
\varphi_{\g} \rangle\times \langle \tau \rangle$ and $\langle
\varphi_{\g} \rangle\times \langle \sigma \rangle$ are normal
subgroups of $\Aut\,E_{\a}$. $\Box$

\begin{co}\label{c3}
The relations (\ref{f11})-(\ref{f13}) are defining relations of
the group $\Aut\,E_{\a}$ with respect to the generators
$\varphi_{\g}$, $\sigma$, and $\tau$ if $\a^3 \neq 1$.
\end{co}

\begin{co}\label{c4}
Let $K$ be an arbitrary field of characteristic $0$ such that the equation $\lambda^2+\lambda+1=0$ has no solution in $K$.
If $\a^3 \neq 1$ (i.e. $\a\neq 1$) then
\bes
\Aut\,E_{\a}\cong K^*\times Z_3.
\ees
\end{co}

When $\a^3 = 1$ the curve $C$ is not elliptic and degenerates into
a product of three linear factors if $\epsilon \in K$ (see (\ref{f3})). Let us
introduce new variables by $u = x+y+\a z$, $v =  \epsilon
x+\epsilon^2 y+\a z$, $w = \epsilon^2 x+\epsilon y+\a z$.
Then $E_{\a}$ is the polynomial algebra $K[u,v,w]$ endowed with the Poisson bracket defined by
\bes \{u, v\} = \mu uv, \,\, \{v, w\} = \mu vw, \,\, \{w, u\} =
\mu wu, \ees where $\mu = 3\a(\epsilon^2 - \epsilon)$. It is easy to check that the maps
 $\psi_{\gamma}$ ($\gamma \in K^*$) and
$\sigma'$ defined by \bes \psi_{\g}(u)=\g u, \ \
\psi_{\g}(v)= v, \ \ \psi_{\g}(w)= w \ees and \bes
\sigma'(u) = v, \ \ \sigma'(v) = w, \ \ \sigma'(w) = u \ees
 are automorphisms of $E_{\a}$. Note that $\sigma=\epsilon^2 \sigma'$.
 Immediate calculations as in the proofs of Theorems \ref{t1} and \ref{t2} give

\begin{theor}\label{t3}
If $\a^3 = 1$ and $\epsilon \in K$ is a root of the equation
$\lambda^2+\lambda+1=0$ then $\psi_{\gamma}$ ($\gamma \in K^*$)
and $\sigma'$  generate $\Aut\,E_{\a}$  and \bes \Aut\,E_{\a}\cong
(K^*)^3\rtimes Z_3. \ees
\end{theor}

If $\a^3 = 1$ and $\epsilon \not\in K$ then $\a=1$. In this case
we cannot use variables $u,v,w$. But we can consider the algebra
$E'_{\a}=E_{\a}\otimes_K K'$ where $K'=K+\epsilon K$ is an
extension of $K$. Every automorphism of $E_{\a}$ can be extended
to an automorphism of $E'_{\a}$ and Theorem \ref{t3} gives the
next result.

\begin{co}\label{c5}
If $\a^3 = 1$ and $\epsilon \not\in K$ then
the group of automorphisms of $E_{\a}$ consist of the
automorphisms given by \bes x \mapsto k_1x + k_2y + k_3z,\ \ y
\mapsto k_1y + k_2z + k_3x, \ \ z \mapsto k_1z + k_2x + k_3y, \ees
where $k_1^3 + k_2^3 + k_3^3 - 3k_1k_2k_3 \neq 0$.
\end{co}
The crucial difference of the case when $\epsilon \not\in K$ is
the absence of the automorphism $\sigma$, just as in the Corollary
\ref{c4}.

\bigskip

\begin{center}
{\bf\large Acknowledgments}
\end{center}

\hspace*{\parindent}

The authors wish to thank several institutions which supported them while they were
working on this project:
Max-Planck Institute
f\"ur Mathematik (the first and the third authors), Department of Mathematics of Wayne
State University in Detroit (the third author), and Instituto de
Matem\'{a}tica e Estat\'{i}stica da Universidade de S\~{a}o Paulo (the first author).


\begin{thebibliography}{99}


\bibitem{Belov-Kontsevich} A.\,Belov-Kanel, M.\,Kontsevich, Automorphisms of the Weyl Algebra,
Letters in Mathematical Physics, 74 (2005), 181--199.



\bibitem{Czer} A.\,J.\,Czerniakiewicz, Automorphisms of a free associative
algebra of rank 2, I, II, Trans. Amer. Math. Soc., 160 (1971),
393--401; 171 (1972), 309--315.


\bibitem{Don} J. Donin and L.\,Makar-Limanov, Quantization of quadratic Poisson brackets on a
polynomial algebra of three variables, J. Pure Appl. Algebra,  129
(1998),  no. 3, 247--261.

\bibitem{DH} J. Dufour et A. Haraki,
 Rotationnels et structures de Poisson quadratiques,  
C. R. Acad. Sci. Paris, V. 312 I (1991), 137-140.

\bibitem{Essen} A.\,van den Essen, Polynomial automorphisms and the
Jacobian conjecture, Progress in Mathematics, 190, Birkhauser verlag,
Basel, 2000.

\bibitem{Freu}   G. Freudenburg, Algebraic theory of locally nilpotent derivations, Springer-Verlag, Berlin, 2006.

\bibitem{Huse} D.\,Husem\"{o}ller, Elliptic Curves, 2nd edition, Springer-Verlag, New York, 2004.

\bibitem{Jung} H.\,W.\,E.\,Jung, \"Uber ganze birationale Transformationen
der Ebene, J. reine angew. Math., 184 (1942), 161--174.


\bibitem{Kulk} W.\,van der Kulk, On polynomial rings in two variables,
Nieuw Archief voor Wiskunde, (3)1 (1953), 33--41.

\bibitem{LX} Z. Liu and P. Xu,
 On quadratic Poisson structures, 
Letters in Math. Phys., 26 (1992), 33-42.


\bibitem{Makar} L.\,Makar-Limanov, The automorphisms of the free algebra with two generators, Funksional.
Anal. i Prilozhen. 4(1970), no.3, 107-108; English translation: in
Functional Anal. Appl. 4 (1970), 262--263.

\bibitem{MakarU2}
L.~Makar-Limanov, U.~U.~Umirbaev, Centralizers in free Poisson algebras,
  Proc. Amer. Math. Soc. 135 (2007), no. 7,  1969--1975.

\bibitem{MLTU}
L.~Makar-Limanov, U.~Turusbekova, U.~Umirbaev,  Automorphisms
and derivations of free Poisson algebras in two variables,  J.
Algebra (submitted).

\bibitem{Nagata} M.\,Nagata, On the automorphism group of $k[x,y]$,
Lect. in Math., Kyoto Univ., Kinokuniya, Tokio, 1972.

\bibitem{Odesskii} A.\,V.\,Odesskii, Elliptic algebras, Russian Math. Surveys, 57 (2002), no. 6, 1127--1162.



\bibitem{Umi24}
I.~P.~Shestakov and U.~U.~Umirbaev, Poisson brackets and two
generated subalgebras of rings of polynomials, Journal of the
American Mathematical Society, 17 (2004), 181--196.

\bibitem{Umi25}
I.~P.~Shestakov and U.~U.~Umirbaev, Tame and wild automorphisms of
rings of polynomials in three variables, Journal of the American
Mathematical Society, 17 (2004), 197--227.

\bibitem{Um19}
 U.~U.~Umirbaev,  The Anick automorphism of free associative algebras,
 J. Reine Angew. Math. 605 (2007), 165--178.


\end{thebibliography}
\end{document}